\magnification=\magstep1
\hsize=16.5 true cm 
\vsize=23.8 true cm
\font\bff=cmbx10 scaled \magstep1
\font\bfff=cmbx10 scaled \magstep2
\font\bffg=cmbx10 scaled \magstep3

\font\smc=cmcsc10 
\parindent0cm
\def\cl{\centerline}           %
\def\bp{\bigskip}              %
\def\mp{\medskip}              %
\def\sp{\smallskip}            %
\def\bc{{\bf c}}               %
           %
\def\Bbb#1{\hbox{\boldmas #1}} %
\def\R{\Bbb R}                 %
\def\N{\Bbb N}                 %
\def\C{\Bbb C}                 %

\expandafter\edef\csname amssym.def
\endcsname{%
       \catcode`\noexpand\@=\the\catcode`\@\space}

\catcode`\@=11

\def\undefine#1{\let#1\undefined}
\def\newsymbol#1#2#3#4#5{\let\next@\relax
 \ifnum#2=\@ne\let\next@\msafam@\else
 \ifnum#2=\tw@\let\next@\msbfam@\fi\fi
 \mathchardef#1="#3\next@#4#5}
\def\mathhexbox@#1#2#3{\relax
 \ifmmode\mathpalette{}{\m@th\mathchar"#1#2#3}%
 \else\leavevmode\hbox{$\m@th\mathchar"#1#2#3$}\fi}
\def\hexnumber@#1{\ifcase#1 0\or 1\or 2\or 3\or 4\or 5\or 6\or 7\or 8\or
 9\or A\or B\or C\or D\or E\or F\fi}

\font\tenmsa=msam10
\font\sevenmsa=msam7
\font\fivemsa=msam5
\newfam\msafam
\textfont\msafam=\tenmsa
\scriptfont\msafam=\sevenmsa
\scriptscriptfont\msafam=\fivemsa
\edef\msafam@{\hexnumber@\msafam}
\mathchardef\dabar@"0\msafam@39
\def\dashrightarrow{\mathrel{\dabar@\dabar@\mathchar"0\msafam@4B}}
\def\dashleftarrow{\mathrel{\mathchar"0\msafam@4C\dabar@\dabar@}}

\def\ulcorner{\delimiter"4\msafam@70\msafam@70 }
\def\urcorner{\delimiter"5\msafam@71\msafam@71 }
\def\llcorner{\delimiter"4\msafam@78\msafam@78 }
\def\lrcorner{\delimiter"5\msafam@79\msafam@79 }
\def\yen{{\mathhexbox@\msafam@55}}
\def\checkmark{{\mathhexbox@\msafam@58}}
\def\circledR{{\mathhexbox@\msafam@72}}
\def\maltese{{\mathhexbox@\msafam@7A}}

\font\tenmsb=msbm10
\font\sevenmsb=msbm7
\font\fivemsb=msbm5
\newfam\msbfam
\textfont\msbfam=\tenmsb
\scriptfont\msbfam=\sevenmsb
\scriptscriptfont\msbfam=\fivemsb
\edef\msbfam@{\hexnumber@\msbfam}
\def\Bbb#1{{\fam\msbfam\relax#1}}
\def\widehat#1{\setbox\z@\hbox{$\m@th#1$}%
 \ifdim\wd\z@>\tw@ em\mathaccent"0\msbfam@5B{#1}%
 \else\mathaccent"0362{#1}\fi}

\def\widetilde#1{\setbox\z@\hbox{$\m@th#1$}%
 \ifdim\wd\z@>\tw@ em\mathaccent"0\msbfam@5D{#1}%
 \else\mathaccent"0365{#1}\fi}
\font\teneufm=eufm10
\font\seveneufm=eufm7
\font\fiveeufm=eufm5
\newfam\eufmfam
\textfont\eufmfam=\teneufm
\scriptfont\eufmfam=\seveneufm
\scriptscriptfont\eufmfam=\fiveeufm

\newsymbol\risingdotseq 133A
\newsymbol\fallingdotseq 133B
\newsymbol\complement 107B
\newsymbol\nmid 232D
\newsymbol\rtimes 226F
\newsymbol\thicksim 2373

\font\eightmsb=msbm8   \font\sixmsb=msbm6   \font\fivemsb=msbm5
\font\eighteufm=eufm8  \font\sixeufm=eufm6  \font\fiveeufm=eufm5
\font\eightrm=cmr8     \font\sixrm=cmr6     \font\fiverm=cmr5
\font\eightbf=cmbx8    \font\sixbf=cmbx6    
      \font\eighti=cmmi8   \font\sixi=cmmi6
\font\ninesy=cmsy9     \font\eightsy=cmsy8  \font\sixsy=cmsy6
     \font\eightit=cmti8  
     \font\eightsl=cmsl8  
     \font\eighttt=cmtt8

\font\eightsmc=cmcsc8
\newskip\ttglue
\newfam\smcfam
\def\eightpoint{\def\rm{\fam0\eightrm}%
  \textfont0=\eightrm \scriptfont0=\sixrm \scriptscriptfont0=\fiverm
  \textfont1=\eighti \scriptfont1=\sixi \scriptscriptfont1=\fivei
  \textfont2=\eightsy \scriptfont2=\sixsy \scriptscriptfont2=\fivesy
  \textfont3=\tenex \scriptfont3=\tenex \scriptscriptfont3=\tenex
  \def\smc{\fam\smcfam\eightsmc}
  \textfont\smcfam=\eightsmc          
\textfont\eufmfam=\eighteufm              \scriptfont\eufmfam=\sixeufm
     \scriptscriptfont\eufmfam=\fiveeufm
\textfont\msbfam=\eightmsb            \scriptfont\msbfam=\sixmsb
     \scriptscriptfont\msbfam=\fivemsb
\def\it{\fam\itfam\eightit}%
  \textfont\itfam=\eightit
  \def\sl{\fam\slfam\eightsl}%
  \textfont\slfam=\eightsl
  \def\bf{\fam\bffam\eightbf}%
  \textfont\bffam=\eightbf \scriptfont\bffam=\sixbf
   \scriptscriptfont\bffam=\fivebf
  \def\tt{\fam\ttfam\eighttt}%
  \textfont\ttfam=\eighttt
  \tt \ttglue=.5em plus.25em minus.15em
  \normalbaselineskip=9pt
  \def\MF{{\manual opqr}\-{\manual stuq}}%
  \let\big=\eightbig
  \setbox\strutbox=\hbox{\vrule height7pt depth2pt width\z@}%
  \normalbaselines\rm}
\def\eightbig#1{{\hbox{$\textfont0=\ninerm\textfont2=\ninesy
  \left#1\vbox to6.5pt{}\right.\n@space$}}}

\catcode`@=13 

\cl{\bffg Proper classes of non-embeddable continua}
\bp
\cl{\bfff Gerald Kuba}
\bp\sp
\vbox{\eightpoint
{\bf Abstract.} Our main result is a construction of a class {\bf H} 
of pathwise connected, locally connected, compact Hausdorff spaces 
such that (i) if $X,Y$ are distinct spaces in {\bf H}
then a continuous, injective mapping from $X$ to $Y$ does not exist;
(ii) {\bf H} contains $2^\kappa$
spaces of weight $\kappa$ and size $\kappa$ 
for every cardinal number $\kappa\geq 2^{\aleph_0}$;
(iii) if $\kappa,\lambda$ are infinite cardinals and $\,\lambda\leq\kappa\,$
then {\bf H} contains $2^\kappa$
spaces of weight $\kappa$ and size $\kappa^\lambda$; 
(iv) the $2^{\aleph_0}$ metrizable continua in {\bf H} are 
subspaces of the plane.
\sp   
}
\bp\sp
{\bff 1. Introduction}
\mp
Write $\,|S|\,$ for the cardinality (the {\it size}) of a set $\,S\,$
and put $\,\bc=|\R|=2^{\aleph_0}\,$.
As usual, $\,w(X)\,$ denotes the {\it weight} of 
the topological space $\,X\,$.
Note that $\,w(X)\leq|X|\leq 2^{w(X)}\,$ 
for every compact Hausdorff space $\,X\,$.
A {\it continuum} 
is a compact and connected Hausdorff space 
with more than one point. Thus $\,|X|\geq \bc\,$ for any continuum $\,X\,$.
If a continuum $\,X\,$ is metrizable 
then $\,|X|=\bc\,$ due to $\,|X|\leq 2^{w(X)}\,$.
The following theorem is classic, see [4].
\mp
(1.1)$\;$ {\it For every infinite cardinal number 
$\kappa$ there exist exactly
$2^\kappa$ continua of weight $\kappa$ up to homeomorphism.}
\mp
The continua constructed in [4] in order to prove (1.1) are not 
pathwise connected. Our first goal is to prove the following theorem
which is also motivated by [3] (see the remark below).
Note that, due to compactness, two continua are {\it non-embeddable} 
if and only if no injective continuous mapping between the spaces
exists. In particular, non-embeddable spaces are non-homeomorphic.
\mp\sp
{\bf Theorem 1.} {\it There exists a 
class $\,{\bf P}\,$ 
of pairwise non-embeddable pathwise connected 
and locally connected continua
such that $\,{\bf P}\,$ contains $2^\kappa$
spaces of weight $\,\kappa\,$ and size $\,\max\{\kappa,\bc\}\,$ 
for every infinite cardinal 
$\,\kappa\,$.}
\mp
It goes without saying that the class $\,{\bf P}\,$ is {\it proper},
that is, $\,{\bf P}\,$ is not a set.
As an immediate consequence of Theorem 1, pathwise connectedness
can be included in (1.1).
An obvious challenge in proving Theorem 1 is to make sure that 
a continuum of arbitratry size $\,\kappa\,$ 
is never embeddable in a continuum of arbitrary size
$\,\lambda>\kappa\,$. 
\mp
An important class of continua, which are not counted up to homeomorphism 
in [4], are the linearly ordered continua.
As usual, a space $\,X\,$ is {\it linearly ordered} 
if the topology of $\,X\,$ is the order topology of some
linear ordering  of $\,X\,$. We point out that 
a continuum is linearly ordered
if and only if it has precisely two non-cut points, see [1] 6.3.8.
It is well known (see [1] 6.3.8.c) that for a linearly ordered 
continuum $\,X\,$
the following three statements are equivalent: (i) $\,X\,$ 
{\it is pathwise connected}, (ii) $\,X\,$ {\it is homeomorphic to} $\,[0,1]\,$,
(iii) $\,w(X)=\aleph_0\,$.
Therefore, in the following theorem the case $\,\kappa=\aleph_0\,$
is excluded.
Note that, of course, every linearly ordered continuum is locally connected.
\mp
{\bf Theorem 2.} {\it There exists a class $\,{\bf L}\,$ 
of pairwise non-embeddable linearly ordered continua
such that $\,{\bf L}\,$ contains $2^\kappa$
spaces of weight $\,\kappa\,$ and size $\,\max\{\kappa,\bc\}\,$ 
for every cardinal $\,\kappa>\aleph_0\,$.}
\mp\sp
In Theorem 2 the restriction $\,\kappa>\aleph_0\,$ does not imply 
$\,\kappa\geq \bc\,$. On the contrary, 
in proving the theorems we have 
to consider the following fact, which is an immediate consequence of 
(9.1) below. 
\mp
(1.2)$\;$ {\it The existence of $\,\bc\,$ uncountable cardinals 
$\,\kappa<\bc\,$ is consistent with {\rm ZFC} set theory.} 
\mp
\eject
{\it Remark.} In [3] we construct a family $\,{\cal G}\,$ of 
pairwise non-embeddable torsion-free groups such that 
$\,{\cal G}\,$ contains $2^\kappa$ groups of cardinality $\kappa$ 
for every infinite cardinal number $\,\kappa\leq\theta\,$ 
where $\,\theta\,$ is  
the first strong limit cardinal of uncountable cofinality.
(Naturally, the numbers $2^\kappa$ are the largest possible.)
A restriction like $\,\kappa\leq\theta\,$ can be avoided 
in Theorem 1.
\bp
{\bff 2. Linearly ordererd continua}
\mp
The clue in proving 
Theorem 1 is to apply the following theorem about totally pathwise disconnected 
linearly ordered continua which is interesting in its own right.
\mp
{\bf Theorem 3.} {\it For every cardinal $\,\kappa>\aleph_0\,$ 
there exist $\,2^\kappa\,$ mutually non-homeomorphic  
continua of weight $\,\kappa\,$ and size 
$\,\max\{\kappa,\bc\}\,$ which are linearly ordered
and totally pathwise disconnected.}
\mp
{\it Remark.} The restriction $\,\kappa>\aleph_0\,$
is inevitable because, as already mentioned,  
every second countable linearly ordered continuum 
is homeomorphic to the compact unit interval $\,[0,1]\,$.
\mp
For the proof of Theorem 3 we need some basic facts 
about linearly ordered continua. First of all we point out 
the following widely known characterization (see [1] 3.12.13 and 6.3.2).
\mp
{\bf Lemma 1.} {\it A space $\,X\,$ is a linearly ordered continuum 
if and only if there exists a linear ordering 
$\,\preceq\,$ of $\,X\,$ such that {\rm (i)} the topology of $\,X\,$
is the order topology induced by $\,\preceq\,$,
{\rm (ii)} for $\,x,y\in X\,$ with $\,x\prec y\,$
there is always a point $\,z\in X\,$ with $\,x\prec z\prec y\,$,
{\rm (iii)} with respect to the ordering $\,\preceq\,$,
every nonempty subset of $\,X\,$
has a supremum and $\,X\,$ has a minimum.}
\mp
Let $\,X\,$ be a linearly ordered continuum. 
If $\,a,b\in X\,$ and $\,a\not= b\,$ 
then the {\it interval $\,I[\{a,b\}]\,$
with endpoints $\,a,b\,$} is the unique subcontinuum $\,Y\,$
of $\,X\,$ such that $\,Y\,$ has precisely two non-cut points and 
these points are $\,a,b\,$. (Alternatively,  $\,I[\{a,b\}]\,$ is 
the intersection of all subcontinua of $\,X\,$ containing 
the set $\,\{a,b\}\,$.) Of course, if $\,S\subset X\,$ and $\,|S|=2\,$ 
then 
\mp
(2.1)\qquad $\;I[S]\,=\,\{\,x\in X\;|\;\min S\preceq x\preceq 
\max S\,\}\;$
\mp
for every linear ordering $\,\preceq\,$ of $\,X\,$ generating 
the topology of $\,X\,$. 
Therefore, {\it in $\,X\,$ a point $x$ lies between points $y$ and $z$}
is a purely topological statement.
Of course, the continuum $\,X\,$ itself is an interval and the two 
endpoints of $\,X\,$ are the only non-cut points of $\,X\,$.
Therefore and in view of (2.1)
there are precisely {\it two} linear orderings of $\,X\,$
which generate the topology of $\,X\,$ and the 
following essential lemma is true.
\mp
{\bf Lemma 2.} {\it If $\,X\,$ and $\,Y\,$ are two linearly ordered continua 
and $\,f\,$ is a bijection from $\,X\,$ onto $\,Y\,$
then $\,f\,$ is a homeomorphism if and only if
$\,f\,$ is either increasing or decreasing
referring to any linear ordering of $\,X\,$ resp.~$\,Y\,$
which generates the topology of $\,X\,$ resp.~$\,Y\,$.}
\mp
We point out that, other than for arbitrary continua, 
the following is true.
\mp
(2.2)\quad {\it The weight of a linearly ordered continuum 
cannot be greater than the size of a dense subset.}
\mp
Indeed, if $\,D\,$ is a dense subset of a linearly ordered continuum $\,X\,$
with endpoints $\,a,b\,$ 
then $\;\{\,I[S]\setminus S\;|\;S\subset D\;\land\;|S|=2\,\}
\cup\{\,I[\{e,x\}]\setminus\{x\}\;|\;e\in\{a,b\}\;\land\;e\not=x\in D\,\}\;$
is a basis of $\,L\,$ of size $\,|D|\,$.
\mp\sp
For the continua proving Theorem 3 we need a special 
linearly ordered continuum as a basic building block.
As usual, $\,\aleph_1\,$ denotes the least cardinal
greater than $\,\aleph_0=|\N|\,$, 
whence $\,\aleph_0<\aleph_1\leq\bc\,$.
It cannot be decided whether
$\,\aleph_1<\bc\,$ or $\,\aleph_1=\bc\,$.
By (1.2) one cannot rule out that 
$\,\aleph_1<\kappa<\bc\,$ for uncountably many cardinals $\,\kappa\,$.
\eject
\mp
{\bf Proposition 1.} {\it There exists
a linearly ordered 
continuum of weight $\,\aleph_1\,$ and size $\,\bc\,$ 
which is first countable
and totally pathwise disconnected.} 
\mp
{\it Remark.} In Proposition 1 the size $\,\bc\,$ is inevitable because 
the size of a first countable compact Hausdorff space 
cannot exceed $\,\bc\,$, see [1] 3.1.30. 
\mp
{\it Proof of Proposition 1.} 
Fix a dense subset $\,D\,$ of $\,[0,1]\,$ 
with $\,\{0,1\}\subset D\,$. Let $\,K[D]\,$ be 
the set of all sequences $\,(x_1,x_2,x_3,...)\,$ of real numbers 
$\,x_n\in[0,1]\,$ such that for every index $\,n\,$ 
the following implication is true.
\sp
(2.3)\qquad\qquad$\;x_n\not\in D\;\;\Longrightarrow\;\;
\forall\,m>n\,:\,\;x_m=0\;$ 
\mp
Consider the subset $\,K[D]\,$ of $\,[0,1]^\N\,$
equipped with the lexicographic ordering $\,\preceq\,$.
So 
\sp
\cl{$\,(x_1,x_2,x_3,...)\prec(y_1,y_2,y_3,...)\,$} 
\sp
for distinct sequences $\,(x_n),(y_n)\,$ 
if and only if $\,x_m<y_m\,$ for 
$\;m\,=\,\min\,\{\,n\in\N\;|\;x_n\not=y_n\,\}\,$.
\mp
We claim that  $\,K[D]\,$ 
is a linearly ordered continuum with respect to the order 
topology induced by $\,\preceq\,$.
Firstly, (ii) in Lemma 1 obviously holds.
Let $\,\pi_n\,$ denote the $n$-th projection
from $\,[0,1]^\N\,$ to $\,[0,1]\,$ and 
let $\,\sup S\,$ denote the supremum of $\,S\subset [0,1]\,$
in the naturally ordered set $\,[0,1]\,$. In particular, 
$\;\sup \emptyset\,=\,0\,$. Then in view of (2.3) it is plain
that if $\,\emptyset\not=Y\subset K[D]\,$ then
the supremum of $\,Y\,$ in $\,K[D]\,$ is the sequence 
$\;(y_1,y_2,y_3,...\,)\;$ where
$\;y_1\,=\,\sup\pi_1(Y)\;$ and 
$\;y_{n+1}\,=\,
\sup\,\pi_{n+1}\big(Y\cap\bigcap_{i=1}^n\pi_i^{-1}(\{y_i\})\big)\;$
for every $\,n\in\N\,$. And, of course, 
$\,(0,0,0,...\,)\,$ is the minimum of $\,K[D]\,$.
\mp
Trivially, $\,|K[D]|=\bc\,$.
Obviously, between two distinct points in $\,K[D]\,$
there always lie $\,|D|\,$ copies of the whole space.
Consequently, every interval in $\,K[D]\,$ is a continuum of weight 
not smaller than $\,|D|\,$. In particular, 
$\,w(K[D])\geq|D|\,$. 
As another consequence, if $\,D\,$ is {\it uncountable} then 
no interval in $\,K[D]\,$
can be homeomorphic to $\,[0,1]\,$ and hence 
the continuum $\,K[D]\,$ is {\it totally pathwise disconnected.}
Considering (2.3), it is evident that 
$\;\bigcup_{n=1}^\infty D^n\!\times\!\{0\}^\N\;$
is a dense subset of $\,K[D]\,$ equipotent with $\,D\,$.
Consequently by (2.2), 
$\,w(K[D])\leq|D|\,$ and hence $\,w(K[D])=|D|\,$.
\mp
The continuum  $\,K[D]\,$ is first countable
because if $\,A\,$ is 
a countable dense subset of $\,D\,$
and $\,x=(x_1,x_2,x_3,...)\,$ is a point 
in $\,K\,$ 
then the countable family of the interiors 
of all intervals with endpoints 
$\;(x_1,...,x_n,a,b,b,b,b,...)\;$
where  $\,n\in\N\,$ is arbitrary and $\,a,b\in\{0,1\}\cup A\,$ 
is a neighborhood basis of $\,x\,$.
Finally, if we choose $\,D\,$ such that $\,|D|=\aleph_1\,$ 
then $\,K[D]\,$ is a continuum exactly as desired,
{\it q.e.d.}
\bp
{\bff 3. Proof of Theorem 3}
\mp
Let $\,\Omega\,$ be the well-ordered
set of all countable ordinal numbers
and $\;{\overline\Omega}\,=\,\Omega\cup\{\omega_1\}\;$
where $\;\omega_1\,=\,\sup\Omega\,=\,\max{\overline\Omega}\;$ 
is the first uncountable ordinal number. (Of course, 
$\,0=\min \Omega=\min {\overline\Omega}\,$.)
Naturally, $\;|{\overline\Omega}|=|\Omega|=\aleph_1\,$.
\mp
Let $\,K\,$ be a linearly ordered continuum of weight $\,\aleph_1\,$
as depicted in Proposition 1.
Let $\,L\,$ be a modification of the classical {\it long line} 
(see [1] 3.12.19)
where the building blocks homeomorphic to $\,[0,1]\,$ 
are replaced by copies of $\,K\,$.
More precisely, put 
\sp
\cl{$\;L\,=\,(\Omega\times(K\setminus\{\max K\})\cup\{(\omega_1,\min K)\}\;$}
\sp
and consider $\,L\,$ equipped with the lexicographic ordering $\,\preceq\,$
with respect to the natural well-ordering of $\,{\overline\Omega}\,$
and the linear ordering of the continuum $\,K\,$ 
which generates the topology of $\,K\,$.
So $\;(\alpha,x)\prec(\beta,y)\;$
if and only if either $\,\alpha<\beta\,$ or $\,\alpha=\beta\,$
and $\,x\in K\,$ is smaller than $\,y\in K\,$.
The maximum of $\,(L,\preceq)\,$ is $\,(\omega_1,\min K)\,$ 
and the minimum of $\,(L,\preceq)\,$ is $\,(0,\min K)\,$.
Similarly as the classical long line, the modified long line 
$\,L\,$ is a linearly ordered continuum of weight $\,\aleph_1\,$ 
and size $\,\bc\,$ and {\it the point $\,(\omega_1,\min K)\,$ 
is the one and only point in $\,L\,$
which does not have a countable neighborhood basis}.
Other than the classical long line, the continuum $\,L\,$ is 
totally pathwise disconnected since its building blocks 
are copies of the totally pathwise disconnected continuum $\,K\,$.
\mp
Let $\,L^*\,$ be a copy of $\,L\,$ 
equipped with the {\it backwards linear ordering} of $\,L\,$.
Formally put $\,L^*=\{-1\}\times L\,$
and declare $\,(-1,x)\in L^*\,$ smaller than $\,(-1,y)\in L^*\,$
if and only if $\,y\in L\,$ is smaller than $\,x\in L\,$.
\mp
Now let $\,\kappa\,$ be an uncountable cardinal number
and 
let $\,\Omega_\kappa\,$ be the well-ordered
set of all ordinal numbers smaller than $\,\kappa\,$
and put $\;{\overline\Omega_\kappa}\,=\,\Omega_\kappa\cup\{\kappa\}\,$.
(Note that $\;\kappa=\max{\overline\Omega_\kappa}=\sup\Omega_\kappa\,$.) 
Let $\,{\cal L}\,$ denote the set of all limit ordinals
in $\,\Omega_\kappa\setminus\{0\}\,$.
Let $\,\Sigma\,$ be the family of all functions from 
$\,{\cal L}\,$ to $\,\{0,1\}\,$.
Since $\,|{\cal L}|=\kappa\,$, we have 
$\;|\Sigma|=2^\kappa\,$.
Our goal is to construct for each $\;\sigma\in\Sigma\;$
a totally pathwise disconnected 
linearly ordered continuum $\,X_\sigma\,$ of weight $\,\kappa\,$
and size $\,\max\{\kappa,\bc\}\,$
such that the spaces $\,X_\sigma\,$ and $\,X_{\sigma'}\,$ are non-homeomorphic
whenever $\;\sigma,\sigma'\in \Sigma\;$ are distinct.
\mp
In order to achieve this, for each $\;\sigma\in\Sigma\;$
we create a sort of hyper-long line $\,X_\sigma\,$
 constructed from 
$\,{\overline\Omega_\kappa}\,$ by placing between each ordinal 
$\,\alpha\in \Omega_\kappa\,$
and its successor $\,\alpha+1\,$ a copy of either the ordered
set $\,L\setminus\{\min L,\max L\}\,$ 
or the ordered set $\,L^*\setminus\{\min L^*,\max L^*\}\,$ 
as follows. With $\;L_0\,=\,L\setminus\{\max L\}\;$ and
$\;L_1\,=\,L^*\setminus\{\max L^*\}\;$ put for each $\,\sigma\in\Sigma\,$
\mp
\cl{$X_\sigma\,=\,\{(\kappa,\min L)\}\,\cup\,
((\Omega_\kappa\setminus{\cal L})\times L_0)\,\cup\,
\bigcup\limits_{\alpha\in {\cal L}}(\{\alpha\}\times L_{\sigma(\alpha)})$}
\sp
and consider $\,X_\sigma\,$ equipped with 
the lexicographic ordering 
with respect to the natural well-ordering of $\,{\overline\Omega_\kappa}\,$
and the linear orderings of $\,L_0\,$ and $\,L_1\,$.
\mp
By Lemma~1, with respect to the order topology,
$\,X_\sigma\,$ is a continuum containing a well-ordered 
copy of $\,\overline{\Omega_\kappa}\,$, whence we assume
$\,\overline{\Omega_\kappa}\subset X_\sigma\,$.
So the minimum $\,(0,\min L)\,$ of $\,X_\sigma\,$
is identified with $\,0=\min\overline{\Omega_\kappa}\,$       
and the maximum $\,(\kappa,\min L)\,$ of $\,X_\sigma\,$
is identified with $\,\kappa=\max\overline{\Omega_\kappa}\,$.     
For every $\,\alpha\in\Omega_\kappa\,$ 
the interval $\,I[\{\alpha,\alpha+1\}]\,$ in $\,X_\sigma\,$
is homeomorphic to the continuum $\,L\,$ and 
order-isomorphic to one of the two linearly ordered sets 
$\,L\,$ and $\,L^*$.
Consequently, $\,X_\sigma\,$ is totally pathwise disconnected
and $\,|X_\sigma|=\max\{\kappa,\bc\}\,$
and $\,w(X_\sigma)=\kappa\,$. 
\mp
The two endpoints $\,0,\kappa\,$ of the 
continuum $\,X_\sigma\,$ can topologically be distinguished. 
For the point $\,0\,$
has a compact neighborhood which is a first countable space, while 
the point $\,\kappa\,$ has no such neighborhood.
(If the cofinality of $\,\kappa\,$ is uncountable 
then there is also the distinction that the point $\,0\,$
has a countable neighborhood basis, while the point $\,\kappa\,$
has not a countable neighborhood basis.)
\mp
Let $\,V_\sigma\,$ be 
the set of all points in the space $\,X_\sigma\,$
which do not have countable neighborhood bases.
Let $\,W_\sigma\,$ be the set of all limit points of $\,V_\sigma\,$
in the space $\,X_\sigma\,$. 
Obviously, $\,W_\sigma\,=\,{\cal L}\,$ 
and hence the set $\,{\cal L}\,$ 
is completely determined by the topology of $\,X_\sigma\,$.  
Moreover, each mapping $\,\sigma\in\Sigma\,$ 
can also be recovered 
from the topology of the continuum $\,X_\sigma\,$.
\mp
Indeed, if $\,\sigma\in\Sigma\,$ and 
$\,\lambda\in{\cal L}\,$ then it is evident that
$\,\sigma(\lambda)=0\,$ if and only if 
the point $\,\lambda\,$ in 
the continuum $\,I[\{\lambda,\kappa\}]\,$ has a countable 
neighborhood basis.
\mp
Consequently, for distinct $\;\sigma,\sigma'\in \Sigma\;$
the continua $\,X_\sigma\,$ and $\,X_{\sigma'}\,$ are 
never homeomorphic and this concludes the proof of Theorem 3.
\eject
\bp
{\bff 4. Proof of Theorem 1} 
\mp
For $\,\kappa\geq\aleph_1\,$ 
let $\,{\cal F}_\kappa\,$ be a family of $\,2^\kappa\,$ 
locally connected continua as depicted in Theorem~3. 
(Of course, every linearly ordered continuum is locally connected.)
Since spaces of distinct weight cannot be homeomorphic, 
the class $\;{\bf F}\,=\,\bigcup\,\{\,{\cal F}_\kappa\;|\;
\kappa\geq \aleph_1\,\}\;$ 
consists of mutually non-homeomorphic 
totally pathwise disconnected linearly ordered continua
and contains $\,2^\kappa\,$ spaces of weight $\,\kappa\,$
and size $\,\max\{\kappa,\bc\}\,$ for every $\,\kappa\geq\aleph_1\,$.
\mp
Our first (and main) goal is to define an operator $\,\Phi\,$ 
on the class $\,{\bf F}\,$
which turns every {\it totally pathwise disconnected} continuum 
$\,X\in{\bf F}\,$ into a {\it pathwise connected} continuum 
$\,\Phi(X)\,$ which is locally connected 
such that $\,|\Phi(X)|=|X|\,$ and $\,w(\Phi(X))=w(X)\,$
and for distinct continua $\,X,Y\in{\bf F}\,$ the 
continua $\,\Phi(X),\Phi(Y)\,$ are non-embeddable.
\mp 
For $\,X\in {\bf F}\,$ consider the compact product space
$\,X\times[-1,1]\,$. 
Let $\,\Psi(X)\,$ be the quotient space
of $\,X\times[-1,1]\,$ modulo the two 
compact subspaces $\,X\times\{-1\}\,$ and $\,X\times\{1\}\,$. 
So $\,X\times\{1\}\,$ resp.~$\,X\times\{-1\}\,$ 
is shrinked to a point $\,p\,$ resp.~$\,q\,$ in $\,\Psi(X)\,$ and 
the subspace $\;X\times{]{-1,1}[}\;$ of $\,X\times[-1,1]\,$
is identical with the subspace 
$\,\Psi(X)\setminus\{p,q\}\,$ of $\,\Psi(X)\,$.
One may picture $\,\Psi(X)\,$ as a {\it double cone} with the {\it basis} 
$\,X\times\{0\}\,$ 
and the {\it upper apex} $\,p\,$ and the {\it lower apex} $\,q\,$.
Putting $\,\ell_x\,=\,\{x\}\times{]{-1,1}[}\,$, 
the {\it rulings} of the double cone $\,\Psi(X)\,$
are the subspaces $\;\ell_x\cup\{p,q\}\;(x\in X)\;$
and they all are homeomorphic copies of $\,[0,1]\,$
connecting the lower apex with the upper apex.
Identifying $\,X\,$ with $\,X\times\{0\}\,$ 
we consider $\,X\,$ as a subset of $\,\Psi(X)\,$.
Of course, $\,\Psi(X)\,$ is a compact Hausdorff space.
Since $\;|X|\geq\bc=|[-1,1]|\,$, we have 
$\;|\Psi(X)|=|X|\,$. In view of [1] 3.1.15,
both apices of the double cone $\,\Psi(X)\,$ have countable 
neighborhood bases. Therefore, since 
$\;w(\Psi(X)\setminus\{p,q\})=w(X\times{]{-1,1}[})=
\max\{w(X),w({]{-1,1}[})\}=w(X)\,$, we have $\;w(\Psi(X))=w(X)\,$.
Furthermore, $\,\Psi(X)\,$ is pathwise connected and locally connected. 
\mp
While {\it lower} and {\it upper} can not be distinguished
topologically, the {\it set} $\,\{p,q\}\,$ is obviously 
characterized as the unique set $\,S\subset\Psi(X)\,$ 
such that $\,|S|=2\,$ and $\,\Psi(X)\setminus S\,$
has infinitely many path-components. 
Of course, the subspace $\;\Psi(X)^*\,=\,\Psi(X)\setminus\{p,q\}\;$
of $\,\Psi(X)\,$ is still connected (but not compact).
The path-components of $\,\Psi(X)^*\,$ 
are the sets $\,\ell_x\;(x\in X)\,$
and they all are homeomorphic to the open unit interval $\,{]0,1[}\,$.  
If $\,a,b\,$ are the two endpoints of the linearly ordered
 continuum $\,X\,$ 
then $\,\Psi(X)^*\setminus\ell_a\,$
and $\,\Psi(X)^*\setminus\ell_b\,$
and $\,\Psi(X)^*\setminus(\ell_a\cup\ell_b)\,$ 
remain connected, whereas $\,\Psi(X)^*\setminus\ell_x\,$ 
is disconnected with two components whenever 
$\,x\in X\setminus\{a,b\}\,$.

\mp
(4.1)\quad {\it The 
topological space $\,X\,$ can be recovered from the topological space
$\,\Psi(X)\,$.} 
\mp
In order to verify (4.1) we firstly 
obtain $\,\Psi(X)^*\,$ from $\,\Psi(X)\,$
by removing the unique set $\,S\,$ of size 2 
where $\,\Psi(X)\setminus S\,$
has infinitely many path-components.  
Then $\,\Psi(X)^*\,$ is homeo\-morphic to $\;X\times{]{0,1}[}\,$.
Since $\,X\,$ is totally pathwise disconnected, 
the path-components of $\,X\times{]{0,1}[}\,$
are the fibers $\;\pi^{-1}(\{x\})\;(x\in X)\;$ where $\,\pi\,$ is the 
canonical projection from $\,X\times{]{0,1}[}\,$ onto $\,X\,$. 
Therefore, the quotient space $\,\Psi(X)^*\,$
{\it modulo pathwise connectedness} must be 
homeomorphic to $\,X\,$ and hence (4.1) is proved.
\mp
So by (4.1) we can be sure that 
for distinct 
continua $\,X,Y\in{\bf F}\,$ the 
continua $\,\Psi(X),\Psi(Y)\,$ are non-homeomorphic. 
(If $\,w(X)\not=w(Y)\,$ then 
$\,\Psi(X),\Psi(Y)\,$ are non-homeomorphic
simply because $\,w(\Psi(X))\not=w(\Psi(Y))\,$.)
However, we can rule out that $\,\Psi(X)\,$ 
is embeddable in $\,\Psi(Y)\,$  only in 
the trivial case $\,w(X)>w(Y)\,$.
To achieve non-embeddability in general we modify  
the double cones $\,\Psi(X)\,$ as follows. 
\mp 
Consider the two disjoint continua $\,[1,2]\,$ and $\,[3,4]\,$ 
in the real line $\,\R\,$.
Let $\,X\in{\bf F}\,$ 
and let $\,a,b\,$ be the two endpoints
of $\,X\,$ and assume that $\,\R\cap\Psi(X)=\emptyset\,$.
Let $\,\Phi(X)\,$ be the topological space obtained from 
$\,X\,$ by attaching the two 
continua $\,[1,2]\,$ and $\,[3,4]\,$ 
to the cone $\,\Psi(X)\,$
at the two points $\,a,b\,$. 
Precisely, $\,\Phi(X)\,$ is the quotient space 
of the topological sum of $\,\Psi(X)\,$ and 
$\,[1,2]\,$ and $\,[3,4]\,$ 
where the point $\,2\,$ is identified with $\,a\,$ 
and the point $\,3\,$ is identified with $\,b\,$.
It is evident that $\,\Phi(X)\,$ is a pathwise connected 
and locally connected compact Hausdorff space 
and $\,|\Phi(X)|=|X|=|\Psi(X)|\,$ and 
$\,w(\Phi(X))=w(X)=w(\Psi(X))\,$.
\mp
The apices $\,p,q\,$ of the double cone $\,\Psi(X)\,$
are points in the space $\,\Phi(X)\,$ 
characterized by the property that $\,\Phi(X)\setminus S\,$
has infinitely many path-components for a set $\,S\,$ with $\,|S|=2\,$
if and only if $\,S=\{p,q\}\,$. And the endpoints $\,a,b\,$
of $\,X\subset\Phi(X)\,$ are 
characterized as the two cut points $\,z_1,z_2\,$ 
of the continuum $\,\Phi(X)\,$
such that no open neighborhood of $\,z_i\,$ 
is homeomorphic to $\,{]0,1[}\,$.
The path-components of $\,\Phi(X)\setminus\{p,q\}\,$
disjoint from $\,\{a,b\}\,$ are the open rulings $\,\ell_x\,$
where $\,x\in X\setminus\{a,b\}\,$. 
There are precisely two path-components which meet $\,\{a,b\}\,$
and these pathwise connected spaces are homeomorphic
to $\;(\{0\}\times{]{-1,1}[})\cup([0,1]\times\{0\})\,$,
a triod minus two non-cut points.
\mp
Now let $\,X_1,X_2\in {\bf F}\,$ and assume that 
$\,f\,$ is a homeomorphism from 
$\,\Phi(X_1)\,$ onto a subspace of $\,\Phi(X_2)\,$.
It is clear that $\,f\,$ must map the 
two apices $\,p_1,q_1\,$ of the double cone $\,\Psi(X_1)\,$ to the two 
apices $\,p_2,q_2\,$ of the double cone $\,\Psi(X_2)\,$ and hence 
$\,f(\Phi(X_1)\setminus\{p_1,q_1\})\,$ is a subset 
of $\,\Phi(X_2)\setminus\{p_2,q_2\}\,$.  
Since $\,f\,$ must map every path-component 
of $\,\Phi(X_1)\setminus\{p_1,q_1\}\,$ to a path-component
of $\,\Phi(X_2)\setminus\{p_2,q_2\}\,$ and since 
$\,\{f(p_1),f(q_1)\}=\{p_2,q_2\}\,$,
we conclude that $\,f(P)\,$ must coincide with precisely one path-component 
of $\,\Phi(X_2)\setminus\{p_2,q_2\}\,$ for every path-component 
$\,P\,$ of $\,\Phi(X_1)\setminus\{p_1,q_1\}\,$
which is not a triod. 
For the two triod path-components $\,T_1,T_2\,$ 
of $\,\Phi(X_1)\setminus\{p_1,q_1\}\,$ it is quite possible 
that $\,f(T_i)\,$ is a proper subset of 
a path-component of $\,\Phi(X_2)\setminus\{p_2,q_2\}\,$
but this path-component must be a triod (minus two non-cut points).
Therefore, since $\,f\,$ has to preserve connectedness,
$\,f(\Phi(X_1))\,$ cannot be disjoint from a path-component 
of $\,\Phi(X_2)\setminus\{p_2,q_2\}\,$. 
Consequently, $\,f(\Psi(X_1))=\Psi(X_2)\,$
because if $\,X\in{\bf F}\,$ and $\,a,b\,$ are the two 
endpoints of $\,X\,$ then $\,\Phi(X)\setminus\{a,b\}\,$ 
has precisely three components 
among which precisely two are homeomorphic with $\,{[0,1[}\,$ 
and the closure of the third one coincides with $\,\Psi(X)\,$.
\mp
So we conclude that the two continua
$\,\Psi(X_1)\,$ and $\,\Psi(X_2)\,$ are homeomorphic 
and hence we obtain $\,X_1=X_2\,$ in view of (4.1).
This finishes the proof of the main part of  Theorem~1, 
namely that 
\sp
\cl{${\bf G}\,=\,\{\,\Phi(X)\;|\;X\in{\bf F}\,\}$}
\sp
is a class consisting of pairwise non-embeddable pathwise connected 
and locally connected continua
such that $\,{\bf G}\,$ contains $2^\kappa$
spaces of weight $\kappa$ and size $\,\max\{\kappa,\bc\}\,$ for every cardinal 
$\,\kappa>\aleph_0\,$.                                     
The proof of Theorem 1 is concluded by adding $\,\bc\,$
appropriate continua of weight $\,\aleph_0\,$ to the class $\,{\bf G}\,$.
This can be accomplished in view of the following theorem.
(Note that every locally connected and {\it metrizable} continuum 
is pathwise connected.)
\mp
{\bf Theorem 4.} {\it There exist $\bc$ pairwise 
non-embeddable locally connected continua 
in the plane $\,\R^2\,$ which are not embeddable in 
any continuum from the class {\bf G}.}
\mp
Theorem 4 can be settled by [5] where $\,\bc\,$ 
pairwise non-embeddable {\it dendrites}
in the plane are constructed which obviously 
are not embeddable in any continuum $\,\Phi(X)\in{\bf G}\,$. 
Nevertheless, in Section 6 we  write down a proof of Theorem 4
which is easier and much shorter than the proof in [5]
and which is also needed in the proof of Theorem~6 below.
\bp
{\bff 5. Proof of Theorem 2}
\mp
A proof of Theorem 2  can easily be extracted from 
the previous proof. Indeed, in view of Lemma 2, 
a moment reflection suffices to see that the class
\sp
\cl{${\bf L}\,=\,\{\,\Phi(X)\setminus(\Psi(X)\setminus X)\;|
\;X\in{\bf F}\,\}$} 
\sp
does the job, {\it q.e.d.}
\bp
{\bff 6. Proof of Theorem 4} 
\mp
In the following, a {\it circle} is a continuum 
homeomorphic to the unit circle $\,x^2+y^2=1\,$ in the plane. 
Obviously,  a continuum $\,X\,$ in the class  $\,{\bf G}\,$ always 
contains infinitely many circles. But such a 
circle must contain the two apices 
of the double cone contained in $\,X\,$.
Thus, no continuum in the class  $\,{\bf G}\,$ 
contains {\it two disjoint} circles.
Therefore we prove 
Theorem~4 by constructing $\,\bc\,$ pairwise 
non-embeddable locally pathwise connected and (path\-wise) connected 
compact subspaces of the plane $\,\R^2\,$ which contain 
infinitely many pairwise disjoint continua homeomorphic to 
the unit circle.
\mp
Define continua $\,Q,A,B\subset\R^2\,$ by 
\sp
\qquad $\,Q\,=\,(\{0,1\}\times[0,1])\,\cup\,
([0,1]\times\{0,{1\over 5},{1\over 4},{1\over 3},{1\over 2},1\})\;$ and
\sp
\qquad $\,A\,=\,(\{0,1\}\times[0,1])\,\cup\,
([0,1]\times\{0,{1\over 4},{1\over 3},{1\over 2},1\})\;$ and
\sp
\qquad $\,B\,=\,(\{0,{1\over 2},1\}\times[0,1])\,\cup\,
([0,1]\times\{0,{1\over 2},1\})\,.$
\mp
In the following, for $\,n\in\N\,$
a {\it $n$-fold} point in a connected space $\,X\,$ 
is a point $\,x\in X\,$ such that every neighborhood of $\,x\,$ 
contains an open neighborhood $\,U\,$ of $\,x\,$ 
where the subspace $\;U\cap X\setminus \{x\}\;$ of $\,X\,$ 
has precisely $\,n\,$ components. 
A point is a {\it multiple point} if it is a $n$-fold point 
for some $\,n\geq 3\,$.
By considering multiple points we can easily verify the following two 
observations. 
\sp
(6.1)\quad {\it the continua $\,A\,$ and $\,B\,$ are non-embeddable.}
\sp
(6.2)\quad {\it the continuum $\,Q\,$ is 
embeddable neither in $\,A\,$ nor in $\,B\,$.}
\mp
Indeed, (6.1) is true because $\,A\,$ has exactly six multiple points 
and all are threefold, whereas 
$\,B\,$ has exactly five multiple points 
among which precisely one is fourfold.
And (6.2) is true because $\,Q\,$ has exactly 
eight multiple points and hence more than the space 
$\,A\,$ and $\,B\,$, respectively.
Obviously, none of the continua $\,Q,A,B\,$ 
has a multiple point which is $n$-fold for $\,n\geq 5\,$.
We also point out the following essential observations.
\mp
(6.3)\quad {\it Each subcontinuum of $\,Q\,$ or of $\,A\,$ 
has either no cut point or $\,\bc\,$ cut points.}
\sp
(6.4)\quad {\it Each subcontinuum of $\,B\,$ 
has either $\,\bc\,$ cut points or at most 
one cut point.
There is precisely one subcontinuum $\,C\,$ of $\,B\,$ 
which has exactly one cut point and contains the point $\,(0,0)$,
namely 
$\;C\,=\,([0,{1\over 2}]^2\setminus{]{0,{1\over 2}}[}^2)\,\cup 
([{1\over 2},1]^2\setminus{]{{1\over 2},1}[}^2)\;$ where
$\,({1\over 2},{1\over 2})\,$ is the only cut point of $\,C\,$.}
\mp
The following statement for $\,C\,$ as above is evident.
\sp
(6.5)\quad {\it Neither $\,Q\,$ nor $\,A\,$ nor $\,B\,$ is embeddable into 
the connected space $\,\bigcup\limits_{n=1}^\infty\big((n,n)+C\big)\,$.}
\sp
Now let $\,{\cal G}\,$ be the family of all mappings 
from $\,\N\,$ to the set $\,\{A,B\}\,$, whence 
$\,|{\cal G}|=2^{\aleph_0}=\bc\,$. 
For each $\,g\in{\cal G}\,$ put 
\sp
\cl{$X_g\;=\;Q\,\cup\,\bigcup\limits_{n=1}^\infty
\big((n,n)+g(n)\big)\,$.}
\sp 
It is evident that $\,X_g\,$ is a 
pathwise connected and locally pathwise connected and closed 
subspace of the plane. Obviously the set of all cut points of 
$\,X_g\,$ is $\,{\cal N}\,=\,\{\,(n,n)\;|\;n\in\N\,\}\,$.
Define $\,{\cal N}^*\,=\,\{\,({1\over 2}+n,{1\over 2}+n)\;|\;n\in\N\,\}\,$.
In view of (6.3) and (6.4) we observe that 
\sp
(6.6)\quad {\it If $\,g\in{\cal G}\,$ and $\,Y\,$ is a 
connected, closed subspace of $\,X_g\,$ with only countably many cut points 
then every cut point of $\,Y\,$ lies in $\,{\cal N}\cup{\cal N}^*\,$.}
\mp
By considering multiple points and 
components of spaces $\,X_g\setminus S\,$ and their closures in $\,X_g\,$
for appropriate sets $\,S\,$ of cut points 
of $\,X_g\,$, in view of (6.1) to (6.6)
it is straightforward to check the following step by step.
\sp
For $\,f,g\in{\cal G}\,$ let $\,h\,$ be a homeomorphism from 
$\,X_f\,$ onto a subspace of $\,X_g\,$. 
Firstly we must have $\,h((1,1))=(1,1)\,$ and $\,h(Q)=Q\,$.
Consequently, $\,h((1,1)+f(1))=(1,1)+g(1)\,$ and $\,h((2,2))=(2,2)\,$
and hence $\,f(1)=g(1)\,$.
Consequently, 
$\,h((2,2)+f(2))=(2,2)+g(2)\,$ and $\,h((3,3))=(3,3)\,$
and hence $\,f(2)=g(2)\,$. And so on, concluding
that $\,f(n)=g(n)\,$ for every $\,n\in\N\,$ or, equivalently, $\,f=g\,$. 
So we conclude that for distinct 
$\,f,g\in{\cal G}\,$ the spaces $\,X_f,X_g\,$ are always 
non-embeddable.
\mp
To gain  compactness, identify the plane $\,\R^2\,$ with 
the field $\,\C\,$ and let $\,\varphi\,$ denote the 
continuous injective mapping from $\,\C\setminus\{-i\}\,$ to $\,\C\,$ 
defined by $\,\varphi(z)=(i+z)^{-1}\,$. 
Then for every $\,g\in{\cal G}\,$ the closure of 
$\,\varphi(X_g)\,$ in $\,\C\,$ is the pathwise connected, 
compact set $\,\{0\}\cup\varphi(X_g)\,$
which clearly is locally connected.
Since $\,0\,$ is the unique point $\,z\,$ in the continuum 
$\,\{0\}\cup\varphi(X_g)\,$ such that 
every neighborhood of $\,z\,$ contains infinitely many cut points
of $\,\{0\}\cup\varphi(X_g)\,$,
the non-embeddability of the spaces 
$\;X_g\;(g\in{\cal G})\;$ implies that 
the $\,\bc\,$ locally connected continua 
$\;\{0\}\cup\varphi(X_g)\;(g\in{\cal G})\;$
are pairwise non-embeddable as well.
Since each continuum contains 
infinitely many pairwise disjoint copies of the unit circle,
the proof of Theorem 4 is finished.
\bp
{\bff 7. Oversized continua}
\mp
As already pointed out, $\,w(X)\leq|X|\,$ 
for every compact Hausdorff space $\,X\,$.
If $\,{\bf P}\,$ and $\,{\bf L}\,$ are classes of continua
as in Theorem 1 and Theorem 2 and $\,X\in{\bf P}\cup{\bf L}\,$ 
and $\,w(X)\geq\bc\,$ then $\,w(X)=|X|\,$. 
So the class $\;\{\,X\in{\bf P}\cup{\bf L}\;|\;w(X)<|X|\,\}\;$ 
is not a proper class but a set of 
size not greater than $\,2^\bc\,$. (The bound $\,2^\bc\,$ 
is trivial and not necessarily sharp. In fact it is sharp 
if and only if $\,2^\kappa=2^\bc\,$ for some cardinal $\,\kappa<\bc\,$.)  
Since the common topologist regards the
weight of a space more important than its size,
a space $\,X\,$ with $\,w(X)<|X|\,$ may be called {\it oversized}.
(If the size of a space is regarded more important than its weight
then oversized continua may be called {\it underweight}.)                   
There arises the question whether a proper class
of pairwise non-embeddable oversized continua exists. 
This question is answered by the following expansion of Theorem 2.
For every infinite cardinal number $\,\kappa\,$ 
let $\,\nu(\kappa)\,$ denote the smallest cardinal $\,\mu\,$
with $\,\kappa^\mu>\kappa\,$. Trivially, $\,\nu(\kappa)\leq\kappa\,$
and hence $\,\kappa<\kappa^{\nu(\kappa)}\leq 2^\kappa\,$.
Of course, if $\,\kappa<\bc\,$ then $\,\nu(\kappa)=\aleph_0\,$
and $\,\kappa^{\nu(\kappa)}=\bc=\max\{\kappa,\bc\}\,$.
\mp
{\bf Theorem 5.} {\it There exists a class $\,{\bf K}\,$ 
of pairwise non-embeddable linearly ordered continua
which contains $2^\kappa$
spaces of weight $\,\kappa\,$ and size $\,\kappa\,$ 
for every cardinal $\,\kappa\geq\bc\,$
and also $2^\kappa$
spaces of weight $\,\kappa\,$ and size $\,\kappa^{\nu(\kappa)}\,$ 
for every cardinal $\,\kappa>\aleph_0\,$.}
\mp
For the proof of Theorem 5 we need the following proposition
whose proof is carried out in Section 8.
\mp
{\bf Proposition 2.} {\it For every cardinal $\,\kappa\geq\bc\,$ 
there exists a linearly ordered continuum $\,L_\kappa\,$ 
such that $\,w(I)=\kappa\,$ 
and $\,|I|=\kappa^{\nu(\kappa)}\,$ for every interval $\,I\,$ 
in $\,L_\kappa\,$.}
\bp
Now, in order to prove Theorem 5, 
let $\,{\bf L}\,$ be the class from Section 5 
which settles Theorem~2. Obviously, we can select 
disjoint subclasses $\,{\bf L}_0,{\bf L}_1\,$
of $\,{\bf L}\,$ such that $\,{\bf L}_i\,$ 
contains $2^\kappa$ continua of weight $\,\kappa\,$ and size 
$\,\max\{\kappa,\bc\}\,$ 
for every cardinal $\,\kappa>\aleph_0\,$ and for both $\,i\in\{0,1\}\,$.
Referring to Proposition 2,
for every cardinal $\,\kappa\geq \bc\,$ 
let $\,L_\kappa\,$ be a 
linearly ordered continuum such that $\,w(I)=\kappa\,$ 
and $\,|I|=\kappa^{\nu(\kappa)}\,$ for every interval $\,I\,$ 
in $\,L_\kappa\,$. 
For every linearly ordered continuum $\,X\,$
let $\,X\lor L_\kappa\,$ denote one of the four continua
which can be created by attaching $\,X\,$ to $\,L_\kappa\,$ such that 
precisely one endpoint of $\,X\,$ is identified 
with precisely one endpoint of $\,L_\kappa\,$. If $\,X\in{\bf L}\,$ 
then both endpoints of $\,X\,$ have 
neighborhoods homeomorphic to $\,[0,1[\,$. 
Therefore, by comparing sizes and weights of intervals, 
the class
\mp
\cl{${\bf K}\,=\,{\bf L}_0\,\cup\,\bigcup\limits_{\kappa\geq\bc}
\{\,X\lor L_\kappa\;|\;X\in{\bf L}_1,\,w(X)=\kappa\,\}$} 
\mp
is one as desired, {\it q.e.d.}
\mp\sp
Since $\,|X|\leq 2^{w(X)}\,$ for every continuum $\,X$,
there arises the question whether
some or every size $\,\kappa^{\nu(\kappa)}\,$ in Theorem 5 
can be replaced with $\,2^\kappa\,$. Of course, 
under the General Continuum Hypothesis 
we have $\,\kappa^{\nu(\kappa)}=2^\kappa\,$ for every infinite cardinal 
$\,\kappa\,$. Furthermore, the class of all cardinals $\,\kappa\,$ with 
$\,\kappa^{\nu(\kappa)}=2^\kappa\,$ is a proper class.
Indeed, by [2] (5.23), if $\,\mu\,$ is an arbitrary cardinal number and 
$\,\mu_1=2^\mu\,$ and $\,\mu_{n+1}=2^{\mu_n}\,$ for every $\,n\in\N\,$  
and if $\,\kappa\,$ is the supremum 
of the set $\;\{\,\mu_n\;|\;n\in\N\,\}\;$ 
then $\,\kappa^{\aleph_0}=2^\kappa\,$ and hence 
$\,\kappa^{\nu(\kappa)}=2^\kappa\,$. 
The situation is more comfortable in the realm of 
the pathwise connected continua. 
In this realm, for every infinite cardinal $\,\kappa\,$ 
there are $\,2^\kappa\,$ non-embeddable continua 
of weight $\,\kappa\,$ and size $\,2^{\kappa}$. 
Moreover, we can expand Theorem 1 as follows. 
(Note that $\,\kappa^1=\kappa\,$ and 
$\,\max\{\kappa^\kappa,\bc\}=2^\kappa\,$ for every infinite 
cardinal $\,\kappa\,$.)
\mp
{\bf Theorem 6.} {\it There exists a class $\,{\bf Q}\,$ 
of pairwise non-embeddable 
pathwise connected and locally connected continua such that
if $\,\kappa\,$ is an infinite cardinal number
and $\,\mu\,$ is a cardinal number with $\,1\leq \mu\leq \kappa\,$
then $\,{\bf Q}\,$ contains $\,2^\kappa\,$
spaces of weight $\,\kappa\,$ and size 
$\,\max\{\kappa^\mu,\bc\}\,$.} 
\mp
To put Theorem 6 into perspective let 
$\,P_\kappa\,$ denote the set of all powers $\,\kappa^\mu\,$
with $\,1\leq \mu\leq \kappa\,$. Then 
$\,|P_\kappa|\geq 2\,$  
since $\,\kappa^1=\kappa<2^\kappa=\kappa^\kappa\,$. 
Of course, $\,\kappa^n=\kappa\,$ for $\,n\in\N\,$
and hence $\;P_\kappa\,=\,
\{\kappa\}\cup\{\,\kappa^\mu\;|\;\aleph_0\leq\mu\leq\kappa\,\}\,$.
Trivially, $\,|P_\kappa|\leq\kappa\,$. 
If $\,\kappa\,$ is a power $\,2^\lambda\,$ then 
$\,\kappa^{\aleph_0}=\kappa\,$.
On the other hand, if $\,\kappa\,$ is of countable cofinality
then $\,\kappa^{\aleph_0}>\kappa\,$ by K\"onig's theorem.
Under the General Continuum Hypothesis
we always have $\;P_\kappa\,=\,\{\kappa,2^\kappa\}\,$.
However, it cannot be ruled out that 
for some proper class {\bf C} of infinite cardinals
we have $\,|P_\theta|=\theta\,$ for every $\,\theta\in{\bf C}\,$.
More precisely, the following is true.
\mp
(7.1) $\;$ {\it It is consistent with {\rm ZFC} that 
if $\,\kappa\geq\aleph_0\,$ is a regular cardinal number
and $\,\theta=2^\kappa\,$ then $\,|P_\theta|=\theta\,$.}
\mp
A short proof of (7.1) is carried out in Section 9.

\bp
In order to prove Theorem 6,  
for every cardinal $\,\kappa\geq\aleph_1\,$ let $\,C_\kappa\,$ 
be a totally pathwise disconnected and locally connected continuum
of weight $\,\kappa\,$ and size $\,\max\{\kappa,\bc\}\,$.
Such continua exist by Theorem 3. 
For $\,\kappa\geq\aleph_1\,$ and any cardinal $\,\lambda\geq 1\,$
consider the topological product $\,(C_\kappa)^\lambda\,$
with precisely $\,\lambda\,$ factors and each factor 
equal to $\,C_\kappa\,$. (If $\,\lambda=1\,$ then   
$\,(C_\kappa)^\lambda=C_\kappa\,$.)
Then $\,(C_\kappa)^\lambda\,$ is a compact Hausdorff space
of weight $\,\max\{\kappa,\lambda\}\,$ 
and size $\,\max\{\kappa^\lambda,\bc\}\,$.
(If $\,\lambda\,$ is infinite then 
$\,\max\{\kappa^\lambda,\bc\}=\kappa^\lambda\,$.) 
Naturally, $\,(C_\kappa)^\lambda\,$ is connected and locally connected 
and totally pathwise disconnected.
Let $\,\Psi_\kappa^\lambda\,$ denote the {\it cone}
with basis $\,(C_\kappa)^\lambda\,$ and apex $\,p\,$ 
obtained from the continuum $\,(C_\kappa)^\lambda\times[0,1]\,$
by shrinking the subcontinuum $\,(C_\kappa)^\lambda\times\{1\}\,$
to a point $\,p\,$. Similarly as for the double cones 
in Section 4, the apex $\,p\,$ has a countable neighborhood basis
in view of [1] 3.1.15. Consequently, since $\,w([0,1])=\aleph_0$,
$\,\Psi_\kappa^\lambda\,$ has weight $\,\max\{\kappa,\lambda\}\,$ 
and size $\,\max\{\kappa^\lambda,\bc\}\,$.
Obviously, the compact Hausdorff space $\,\Psi_\kappa^\lambda\,$
is pathwise connected and locally connected
and has no cut points. 
The basis $\,(C_\kappa)^\lambda\times\{0\}\,$ 
of the cone $\,\Psi_\kappa^\lambda\,$ is 
totally pathwise disconnected.   
Consequently, if $\,p\,$ is the apex of
the cone $\,\Psi_\kappa^\lambda\,$ then the path-components of 
the connected subspace $\,\Psi_\kappa^\lambda\setminus\{p\}\,$
of $\,\Psi_\kappa^\lambda\,$ are the sets 
$\;\{x\}\times{[0,1[}\;(x\in C_\kappa)\;$ and hence 
they are homeomorphic to $\,{[0,1[}\,$.
Furthermore, if $\,p\not=y\in\Psi_\kappa^\lambda\,$
 then $\,\Psi_\kappa^\lambda\setminus\{y\}\,$
has either one or two path-components. 
So $\,\Psi_\kappa^\lambda\setminus\{z\}\,$ has infinitely 
many path-components if and only if $\,z=p\,$.
\mp
Now let $\,{\bf G}\,$ be the class of continua 
defined in Section 4. Each $\,X\in{\bf G}\,$ 
is a double cone with two copies of $\,[0,1]\,$
attached. (Note that $\,X\,$ contains infinitely many 
circles, but two circles are never disjoint.)
There are precisely two 
non-cut points $\,z\,$ in $\,X\,$ such that 
some open neighborhood of $\,z\,$ is homeomorphic to $\,{[0,1[}\,$. 
Let $\,a_X\,$ be one of these two points in $\,X\,$ 
for every $\,X\in{\bf G}\,$. 
Now if $\,X\in {\bf G}\,$ and $\,\Psi_\kappa^\lambda\,$ 
is a cone as above (and assumed to be disjoint from $\,X\,$)
then let $\,\sigma(X,\Psi_\kappa^\lambda)\,$
denote the continuum obtained by 
attaching $\,\Psi_\kappa^\lambda\,$ to $\,X\,$ 
such that the special point $\,a_X\,$ is identified with the apex of the cone 
$\,\Psi_\kappa^\lambda\,$.
Naturally, $\,\sigma(X,\Psi_\kappa^\lambda)\,$
is a pathwise connected and locally connected continuum
of weight $\,\max\{w(X),w(\Psi_\kappa^\lambda)\}\,$ 
and size $\,\max\{|X|,|\Psi_\kappa^\lambda|\}\,$.
So if $\,w(X)=\kappa\,$ (and hence $\,|X|=\max\{\kappa,\bc\}\,$) 
then $\,w(\sigma(X,\Psi_\kappa^\lambda))=\kappa\,$ 
for every $\,\lambda\leq\kappa\,$
and $\,|\sigma(X,\Psi_\kappa^\lambda)|=\kappa^\lambda\,$ 
for every infinite $\,\lambda\leq\kappa\,$
and $\,|\sigma(X,\Psi_\kappa^n)|=\max\{\kappa,\bc\}\,$ 
for every $\,n\in\N\,$.
\mp
With $\,\nu(\kappa)\,$ defined as above, 
we have $\;P_\kappa\,=\,\{\kappa\}\cup
\{\,\kappa^\mu\;|\;\nu(\kappa)\leq\mu\leq\kappa\,\}\;$
for every infinite cardinal $\,\kappa\,$. 
For every $\,\theta\in P_\kappa\,$
let $\,\mu(\theta)\,$ be the smallest cardinal $\,\lambda\,$
with $\,\kappa^\lambda=\theta\,$. So $\,\mu(\theta)=1\,$ 
for $\,\theta=\kappa\,$
and $\,\nu(\kappa)\leq\mu(\theta)\leq\kappa\,$ for 
$\,\theta\in P_\kappa\setminus\{\kappa\}\,$
and $\;|\{\,\mu(\theta)\;|\;\theta\in P_\kappa\,\}|=|P_\kappa|\,$.
\mp
For every $\,\kappa\geq\aleph_1\,$ 
the class $\;{\cal G}_\kappa\,=\,\{\,X\in{\bf G}\;|\;w(X)=\kappa\,\}\;$
is a set of size $\,2^\kappa\,$. Since 
$\,|P_\kappa|\leq\kappa\leq|{\cal G}_\kappa|\,$
and $\,{\cal G}_\kappa\cap{\cal G}_{\kappa'}=\emptyset\,$  whenever 
$\,\kappa\not=\kappa'$, 
for every pair $\,(\kappa,\theta)\,$ 
with $\,\kappa\geq\aleph_1\,$ and $\,\theta\in P_\kappa\,$ 
we can choose a subfamily $\,{\cal G}_\kappa^\theta\,$
of $\,{\cal G}_\kappa\,$ such that 
$\,|{\cal G}_\kappa^\theta|=2^\kappa\,$
and $\,{\cal G}_\kappa^\theta\cap{\cal G}_{\kappa'}^{\theta'}=\emptyset\,$
whenever $\,(\kappa,\theta)\not=(\kappa',\theta')\,$. 
\mp
Now let $\,X,Y\in{\bf G}\,$ 
and let $\,\Psi_\kappa^\lambda,\Psi_{\kappa'}^{\lambda'}\,$ 
be two cones and assume that $\,f\,$ is a homeomorphism 
from $\,\sigma(X,\Psi_\kappa^\lambda)\,$ to a subcontinuum 
of $\,\sigma(Y,\Psi_{\kappa'}^{\lambda'})\,$. 
The point $\,a_X\,$ in the continuum 
$\,\sigma(X,\Psi_\kappa^\lambda)\,$
is the unique point $\,z\,$ with the property that 
the whole space minus the singleton $\,\{z\}\,$ 
splits into precisely two components, one component is pathwise connected
{\it and contains a circle},
the other splits into infinitely many path-components 
homeomorphic to $\,{[0,1[}\,$.
If any subcontinuum $\,S\,$ of $\,\sigma(Y,\Psi_{\kappa'}^{\lambda'})\,$
contains such a special point $\,z\,$ 
(and hence a circle)
then there is obviously no other possibility than $\,z=a_Y\,$.
Consequently, $\,f(a_X)=a_Y\,$
and $\,f\,$ maps every path-component 
of $\,\sigma(X,\Psi_\kappa^\lambda)\setminus\{a_X\}\,$ to a 
path-component of $\,\sigma(Y,\Psi_{\kappa'}^{\lambda'})\setminus\{a_Y\}\,$. 
Of course, the path-component $\,X\setminus\{a_X\}\,$ 
is not embeddable in $\,{[0,1[}\,$ and hence 
$\;f(X\setminus\{a_X\})\subset Y\setminus\{a_Y\}\;$
and hence $\,f(X)\subset Y\,$ and hence $\,X=Y\,$. 
Therefore we can be sure that if $\,X,Y\,$
are  distinct continua 
in the class $\,{\bf G}\,$  
and $\,\Psi_\kappa^\lambda,\Psi_{\kappa'}^{\lambda'}\,$ 
are arbitrary cones 
then the continua 
$\,\sigma(X,\Psi_\kappa^\lambda),\sigma(Y,\Psi_{\kappa'}^{\lambda'})\,$ 
are non-embeddable. 
Consequently, for every cardinal $\,\kappa\geq\aleph_1\,$
and every $\,\theta\in P_\kappa\,$ the family 
\sp
\cl{$\,{\cal Y}(\kappa,\theta)\,:=\,
\{\,\sigma(X,\Psi_\kappa^{\mu(\theta)})\;|\;X\in{\cal G}_\kappa^\theta\,\}$}
\sp
consists of $\,2^\kappa\,$ pathwise connected and locally connected continua
of weight $\,\kappa\,$ and size $\,\theta\,$, 
and every union of families $\,{\cal Y}(\cdot,\cdot)\,$
consists of pairwise non-embeddable continua.
It is clear that no continuum containing at least two {\it disjoint} 
circles 
is embeddable into any continuum in any family $\,{\cal Y}(\cdot,\cdot)\,$.
(Every continuum in $\,{\cal Y}(\cdot,\cdot)\,$
contains infinitely many circles but two circles are never disjoint.)
Therefore, if we define $\,{\bf H}\,$ as the union of 
\sp
\cl{$\;\bigcup\,\{\,{\cal Y}(\kappa,\theta)\;|\;\kappa\geq\aleph_1\;
\land\;\theta\in P_\kappa\,\}\;$}
\sp
and the family of all continua in the plane 
which prove Theorem 4 in Section 6 then the proof of Theorem 6 
is concluded.
\bp\mp
{\bff 8. Proof of Proposition 2}
\mp
For a limit ordinal $\,\lambda>0\,$ 
(which may be an infinite cardinal number)
and a linearly ordered set $\,L\,$
consider the set $\,L^\lambda\,$
of all $\lambda$-sequences in $\,L\,$   
equipped with  the lexicographic ordering.
So $\,a\in L^\lambda\,$
when $\;a\,=\,(x_\alpha)_{\alpha<\lambda}\;$ is a mapping from 
the set of all ordinals smaller than $\,\lambda\,$ to 
$\,L\,$ and $\;a=(x_\alpha)_{\alpha<\lambda}\;$ is smaller 
than $\;b=(y_\alpha)_{\alpha<\lambda}\;$ if and only if $\,a\not=b\,$
and $\,x_\beta\,$ is smaller than $\,y_\beta\,$ for the least ordinal
$\,\beta\,$ where $\,x_\beta\not=y_\beta\,$.
\sp
Now assume that $\,L\,$ is a linearly ordered continuum of 
size and weight $\,\kappa\,$. (Such a continuum exists 
for many reasons, for example by Theorem 2.)
Put $\,\lambda=\nu(\kappa)\,$. 
The proof is finished by verifying that 
the definition $\;L_\kappa\,:=\,L^\lambda\;$ does the job.
First of all, the size of $\,L_\kappa\,$
equals $\,|L|^\lambda=\kappa^{\nu(\kappa)}\,$. 
\sp
Clearly, the constant sequence $\,(\max L)\,$ resp.~$\,(\min L)\,$
is the maximum resp.~the minimum of $\,L_\kappa\,$.
Let $\,\sup A\,$ denote the supremum of $\,A\subset L\,$
in $\,L\,$ with the convention 
$\;\sup \emptyset=\min L\,$. For every ordinal $\,\alpha<\lambda\,$
let $\;\pi_\alpha\;$ denote the projection
$\;(x_\beta)_{\beta<\lambda}\mapsto x_\alpha\;$
from $\,L_\kappa\,$ onto $\,L\,$.
Similarly as in the proof of Proposition 1,
for $\;\emptyset\not= Y\subset L_\kappa\;$
define $\;y\,=\,(y_\beta)_{\beta<\lambda}\;$ recursively via
$\;y_0\,:=\,\sup\pi_0(Y)\;$ 
and 
\sp
\cl{$\;y_{\beta}\,:=\,\sup\,\pi_{\beta}
\big(Y\cap\bigcap\limits_{\alpha<\beta}\pi_\alpha^{-1}(\{y_\alpha\})\big)\;$}
\sp
for every ordinal number $\,\beta\,$ with $\,0<\beta<\lambda\,$. 
A moment's reflection suffices to see that $\,y\,$
is the supremum of $\,Y\,$.
\sp
Obviously, between any two points in $\,L_\kappa\,$
there lie at least $\,|L|\,$ mutually disjoint 
order-isomorphic (and homeomorphic) copies
of $\,L_\kappa\,$. In particular, $\,L_\kappa\,$ has no pair of consecutive 
points and hence, by Lemma 1, $\,L_\kappa\,$
is a linearly ordered continuum. As a second consequence,
the weight of $\,L_\kappa\,$ cannot be smaller than $\,|L|=\kappa\,$.
\sp
Let $\,E\,$ be the set of all elements of $\,L_\kappa\,$
which are {\it eventually constant} $\lambda$-sequences. 
So $\,(x_\alpha)_{\alpha<\lambda}\,$ lies in $\,E\,$
if and only if for some ordinal $\,\beta<\lambda\,$ and some point 
$\,a\in L\,$ we have $\;x_\alpha=a\;$ whenever
$\;\beta\leq\alpha<\lambda\,$.
Clearly $\,E\,$ is a dense subset of the continuum $\,L_\kappa\,$.
Since there are precisely $\,|L|=\kappa\,$ choices for $\,a\in L\,$
and since $\,\kappa^\mu=\kappa\,$ 
for every cardinal $\,\mu<\lambda=\nu(\kappa)\,$
(and since $\,\lambda\,$ is the total 
number of ordinals $\,\beta<\lambda\,$),  
we have $\;|E|=\kappa\,$. So the weight of $\,L_\kappa\,$
cannot be greater than $\,\kappa\,$ by (2.2).
Thus we derive $\,w(L_\kappa)=\kappa\,$. Moreover, since every interval in 
$\,L_\kappa\,$ contains a copy of $\,L_\kappa\,$, 
we have $\,w(I)=\kappa\,$ and 
$\,|I|=|L_\kappa|=\kappa^{\nu(\kappa)}\,$
for every interval $\,I\,$ in $\,L_\kappa\,$.
This concludes the proof of Proposition 2.
\bp
{\bff 9. Proof of (7.1)}
\mp
Obviously, (7.1) is an immediate consequence of the following 
consistency result which is interesting in its own right.
\mp
(9.1)$\;$ {\it It is consistent with {\rm ZFC} that 
the continuum function $\,\kappa\mapsto 2^\kappa\,$
is strictly increasing and that for every regular 
cardinal $\,\kappa\geq\aleph_0\,$ there exist precisely
$\,2^\kappa\,$ cardinals $\,\lambda\,$ with
$\;\kappa<\lambda<2^\kappa\,$.}
\mp
{\it Proof.} With $\,\kappa,\lambda,\mu\,$ denoting infinite cardinals,
define in G\"odel's universe L for every regular $\,\kappa\,$
a cardinal number $\,\theta(\kappa)\,$ by 
$\;\theta(\kappa)\,:=\,\min\,\{\,\mu\;|\;\mu=\aleph_\mu \;\land\;
{\rm cf}\,\mu=\kappa^+\,\}\,$.
Then $\;|\{\,\lambda\;|\;\kappa<\lambda<\theta(\kappa)\,\}|=\theta(\kappa)\;$ 
holds in every generic extension of L.
By applying Easton's theorem [1, 15.18] one can 
create an Easton universe E generically extending L
such that the continuum function $\;\kappa\mapsto 2^\kappa=\kappa^+\;$
in L is changed into $\;\kappa\mapsto 2^{\kappa}=g(\kappa)\;$ in E
with $\,g(\kappa)=\theta(\kappa)\,$ for every regular $\,\kappa\,$.
So in E we have $\;|\{\,\lambda\;|\;\kappa<\lambda<2^\kappa\,\}|=2^\kappa\;$ 
for every regular $\,\kappa\,$. By definition, in E we have 
$\,2^\kappa<2^\lambda\,$ whenever $\,\kappa,\lambda\,$ are regular 
with $\,\kappa<\lambda\,$. As in every Easton universe
(see [2] Exercise 15.12),
in E the {\it Singular Cardinal Hypothesis} 
holds, that is, if 
$\,2^{{\rm cf}\kappa}<\kappa\,$ then $\,\kappa^{{\rm cf}\kappa}=\kappa^+\,$.
Therefore, in view of [2] Theorem 5.22, 
if $\,\mu\,$ is singular in E
then $\,2^\mu=(2^{<\mu})^+\,$ is a successor cardinal in E
(because the continuum function is already forced  
to be strictly increasing on the {\it regular} cardinals),
whereas $\,2^\kappa\,$ is {\it by definition} 
a limit cardinal in E for every regular $\,\kappa\,$
in E. Consequently, the continuum function is strictly 
increasing on all cardinals in E, {\it q.e.d.}

\bp\bp\bp
{\bff References}
\mp
[1] R.~Engelking,  {\it General Topology, revised and completed edition.}
Heldermann 1989. 
\sp
[2] T.~Jech, {\it Set Theory}, 3rd ed. Springer 2002.
\sp
[3] G.~Kuba, {\it Many non-embeddable groups.} 
J.~Group Theory {\bf 29} (2026), 933-945. 
\sp
[4] F.W.~Lozier and R.H.~Marty, {\it The number of continua.}
Proc.~American Math.Soc.~{\bf 40} 

\rightline{(1973) 271-273.}
\sp
[5] K.~Sieklucki, {\it On a family of power c consisting of R-uncomparable 
dendrites.} 

\rightline{Fund.~Math. {\bf 46} (1959), 330-335.}

\bp\bp
Gerald Kuba

Institute of Mathematics, BOKU University, Vienna

{\sl E-mail:} {\tt gerald.kuba(at)boku.ac.at}
\end